%% file: main.tex
\documentclass{amsart}

\title{Did Turing prove the undecidability of the halting problem?}

\author[Hamkins]{Joel David Hamkins}
\address[Joel David Hamkins]
{O'Hara Professor of Logic, Department of Philosophy, University of Notre Dame, 100 Malloy Hall, Notre Dame, IN 46556 USA 
}
\email{jdhamkins@nd.edu}
\urladdr{http://jdh.hamkins.org}

\author[Nenu]{Theodor Nenu}
 \address[Theodor Nenu]
         {Early Career Research Fellow in AI and Theoretical Philosophy, Institute for Ethics in AI, Faculty of Philosophy, Balliol College, University of Oxford, Oxford, United Kingdom}
\email{theodor.nenu@philosophy.ox.ac.uk}
\urladdr{https://www.philosophy.ox.ac.uk/people/theodor-nenu}

\thanks{Commentary can be made about this article on the first author's blog at \href{https://jdh.hamkins.org/turing-halting-problem}{https://jdh.hamkins.org/turing-halting-problem}. The second author would like to thank Peter Millican for many fruitful conversations on Alan Turing's work over the years, and for sparking his initial interest in the titular question.}

\usepackage{latexsym,amsfonts,amsmath,amssymb,mathrsfs}
\usepackage{changepage}
\usepackage{csquotes}
\usepackage{braket}
\usepackage{bbm}
\usepackage[hidelinks]{hyperref}
\usepackage{relsize}
\usepackage[cmyk,svgnames,dvipsnames]{xcolor}
\usepackage{tikz}
\usetikzlibrary {3d}
\usepackage{comment}
\usetikzlibrary{arrows,arrows.meta,petri,topaths,positioning,shapes,shapes.misc,patterns,calc,decorations.pathreplacing,hobby}
\pgfdeclarelayer{foreground}
\pgfdeclarelayer{background}
\pgfdeclarelayer{boardarrows}
\pgfdeclarelayer{boardgrid}
\pgfdeclarelayer{boardshades}
\pgfsetlayers{boardshades,boardgrid,boardarrows,background,main,foreground}
\usepackage{wrapfig} 
\usepackage{float}
\usepackage[utf8]{inputenc}
\RequirePackage{doi}

\usepackage{enumitem}

\include{MathMacrosJDH}
\usepackage[backend=bibtex,style=alphabetic,maxbibnames=15,maxcitenames=6,dateabbrev=false]{biblatex}
\renewcommand{\UrlFont}{} 
\renewbibmacro{in:}{\ifentrytype{article}{}{\printtext{\bibstring{in}\intitlepunct}}} 
\DeclareFieldFormat{url}{\UrlFont\url{#1}} 
\DeclareFieldFormat{urldate}{
  (version \thefield{urlday}\addspace%
  \mkbibmonth{\thefield{urlmonth}}\addspace%
  \thefield{urlyear}\isdot)}
\DeclareFieldFormat{eprint:arxiv}{
  \ifhyperref
    {\href{http://arxiv.org/abs/#1}{%
        arXiv\addcolon\nolinkurl{#1}}\iffieldundef{eprintclass}{}{\UrlFont{\mkbibbrackets{\thefield{eprintclass}}}}}
    {arXiv\addcolon\nolinkurl{#1}\iffieldundef{eprintclass}{}{\UrlFont{\mkbibbrackets{\thefield{eprintclass}}}}}}

\begin{document}

\begin{abstract}

We discuss the accuracy of the attribution commonly given to Turing \cite{Turing1936:On-computable-numbers} for the computable undecidability of the halting problem, coming eventually to a nuanced conclusion.

\end{abstract}
\maketitle

\begin{quote}
\setcounter{tocdepth}{1}
\footnotesize
\tableofcontents
\end{quote}

The \emph{halting problem} is the decision problem of determining whether a given computer program halts on a given input, a problem famously known to be computably undecidable. In the computability theory literature, one quite commonly finds attribution for this result given to Alan Turing \cite{Turing1936:On-computable-numbers}, and we should like to consider the extent to which these attributions are accurate. After all, the term \emph{halting problem}, the modern formulation of the problem, as well as the common self-referential proof of its undecidability, are all---strictly speaking---absent from Turing's work. Indeed, Turing's machines as he presents them are not designed ever to halt, but rather specifically to run forever. However, Turing does introduce the concept of an undecidable decision problem, essentially proving that the problems of deciding whether a given program is circle-free or whether it will ever print a specific symbol are both undecidable: we will henceforth call these decision problems the \textit{circle-free problem} and the \textit{printing problem}. This latter problem is easily seen to be computably equivalent to the halting problem and can arguably serve in diverse contexts and applications in place of the halting problem---they are easily translated to one another. Furthermore, Turing laid down an extensive framework of ideas sufficient for the contemporary analysis of the halting problem, including: the definition of Turing machines; the labeling of programs by numbers in a way that enables programs to be enumerated and also for them to be given as input to other programs; the existence of a universal computer; the undecidability of several problems that, like the halting problem, take other programs as input, including the circle-free problem, the printing problem, and the infinite-printing problem, as well as the Hilbert-Ackermann Entscheidungsproblem. In light of these facts, and considering some general cultural observations, by which mathematical attributions are often made not strictly for the exact content of original work, but also generously in many cases for the further aggregative insights to which those ideas directly gave rise, ultimately we do not find it unreasonable to offer qualified attribution to Turing for the undecidability of the halting problem. That said, we would also find it incorrect for an author to suggest that one will find a discussion of the halting problem or a proof of its undecidability in \cite{Turing1936:On-computable-numbers} or that one will find there anything like the familiar self-referential argument for the undecidability of the halting problem.

\section{The halting problem}\label{Section.Halting-problem}

The halting problem can be formulated for essentially any desired notion of computability---Turing machines, multi-tape big-alphabet Turing machines, register machines, flowchart machines, oracle machines, or what have you---and in each case, the main observation to be made is that there is no computable procedure using that notion of computability to determine whether a given program of that sort will halt on a given input.

One of the most commonly seen contemporary arguments proceeds as follows. We suppose toward contradiction that there is a computational procedure to decide whether or not a given program halts on a given input. We fix a particular such procedure, which we shall use as a subroutine in the following algorithm $q$. What program $q$ does on input $p$, a program, is that it asks the subroutine whether $p$ would halt if given $p$ itself as input. If yes, then our program $q$ goes immediately into a nonhalting infinite loop; if no, then our program $q$ halts immediately. In each case, therefore, program $q$ performs the opposite halting behavior on input $p$ than $p$ itself would. It follows from this that if we should run program $q$ on input $q$, then we would see that $q$ halts on $q$ if and only if $q$ does not halt on $q$, which is a contradiction. So there can be no computational procedure to decide the halting problem. 

\section{Did Turing prove the undecidability of the halting problem?}\label{Section.Quotes-giving-Turing-attribution}

One quite commonly finds in the research literature the attribution for the undecidabilty of the halting problem given to Turing and specifically to his 1936 paper, \enquote{On computable numbers with an application to the Entscheidungsproblem} \cite{Turing1936:On-computable-numbers}. Let us provide a number of suggestive quotations in support of this claim.

\medskip
\begin{enumerate}\parskip=8pt\small

    \item John Stillwell (\cite[p. 370]{stillwell2022story}, bold original): \enquote{The problem of deciding whether a given machine halts of a given input---the so-called \textbf{halting problem}---must be unsolvable. This result was also observed by Turing (1936).}

\enlargethispage{20pt}
     \item Graham Priest (\cite[pp. 105-107]{priest2017logic}): \enquote{Is there an algorithm we can apply to a program (or, more precisely, its code number) and inputs, to determine whether or not a computation with that program and those inputs terminates? The answer is \textit{no}. And this is what Turing proved. (...) The result is known as the \textit{Halting Theorem}}\goodbreak

      \item Scott Aaronson (\cite[p. 21]{aaronson2013quantum}, \cite{aaronson1999can}): 
      \begin{enumerate}\parskip=4pt
          \item \enquote{Turing's first result is the existence of a \enquote{universal} machine (...) But this result is not even the main result of the paper. So what is the main result? It's that there's a basic problem, called the halting problem, that no program can ever solve.}

          \item \enquote{Turing proved that this problem, called the Halting Problem, is unsolvable by Turing machines. The proof is a beautiful example of self-reference. It formalizes an old argument about why you can never have perfect introspection: because if you could, then you could determine what you were going to do ten seconds from now, and then do something else. Turing imagined that there was a special machine that could solve the Halting Problem. Then he showed how we could have this machine analyze itself, in such a way that it has to halt if it runs forever, and run forever if it halts.} 
      \end{enumerate}

    \item Thomas Cormen (\cite[p. 210]{cormen2013algorithms}, bold original): \enquote{(T)here are problems for which it is provably impossible to create an algorithm that always gives a correct answer. We call such problems \textit{\textbf{undecidable}}, and the best-known one is the \textit{\textbf{halting problem}}, proven undecidable by the mathematician Alan Turing in 1937. In the halting problem, the input is a computer program A and the input $x$ to A. The goal is to determine whether program A, running on input $x$, ever halts.}

    \item Dexter Kozen (\cite[pp. 243-244]{kozen2012automata}): \enquote{The technique of diagonalization was first used by Cantor to show that there are fewer real algebraic numbers than real numbers. Universal Turing Machines and the application of Cantor's diagonalization technique to prove the undecidability of the halting problem appear in Turing's original paper.} 

    \item Oron Shagrir (\cite[p. 3]{shagrir2006godel}): \enquote{[In his 1936 paper] Turing provides a mathematical characterization of his machines, proves that the set of these machines is enumerable, shows that there is a universal (Turing) machine, and describes it in detail. He formulates the halting problem, and proves that it cannot be decided by a Turing machine. On the basis of that proof, Turing arrives, in section 11, at his ultimate goal: proving that the \textit{Entscheidungsproblem} is unsolvable.}

    \item Douglas Hofstadter (\cite[p. XII]{hofstadter2004turing}): \enquote{Fully to fathom even one other human being is far beyond our intellectual capacity --- indeed, fully to fathom even one's own self is an idea that quickly leads to absurdities and paradoxes. This fact Alan Turing understood more deeply than nearly anyone ever has, for it constitutes the crux of his work on the halting problem.}

    \item Roger Penrose (\cite[p. 30, our emphasis]{penrose1994shadows}): \enquote{The mathematical proofs that Hilbert's tenth problem and the tiling problem are not soluble by computational means are difficult, and I shall certainly not attempt to give the arguments here. The central point of each argument is to show, in effect, how any Turing-machine action can be coded into a Diophantine or tiling problem. \textit{This reduces the issue to one that Turing actually addressed in his original discussion: the computational insolubility of the halting problem}.}

\enlargethispage{20pt}%
    \item Piergiorgio Odifreddi (\cite[p. 150]{odifreddi1992classical}):
    \enquote{The name [of the following theorem] comes from its original formulation, which was in terms of Turing machines, and in that setting it shows that there is no Turing machine that decides whether a universal Turing machine halts or not on given arguments.\\[1ex] 
    \textbf{Theorem II.2.7 Unsolvability of the Halting Problem (Turing [1936])} \textit{The set defined by $\braket{x,e} \in \mathcal{K}_0 \iff x \in \mathcal{W}_e \iff \varphi_e(x)\downarrow$ is r.e. and nonrecursive.}}

    \item Hartley Rogers, Jr. (\cite[p. 19]{rogers1987theory}, underline original): \enquote{There is no effective procedure by which we can tell whether or not a given effective computation will eventually come to a stop. (Turing refers to this as the unsolvability of the \underline{halting problem} for machines. This and the existence of the universal machine are the principal results of Turing's first paper.)}
\medskip
\end{enumerate}
Finally, the first author of this paper:
\medskip
\begin{enumerate}[resume]\small
    \item  Joel David Hamkins (\cite[\S6.5]{Hamkins2021:Lectures-on-the-philosophy-of-mathematics}): ``Is the halting problem computably decidable? In other words, is there a computable procedure, which on input $(p,n)$ will output yes or no depending on whether program $p$ halts on input $n$? The answer is no, there is no such computable procedure; the halting problem for Turing machines is not computably decidable. This was proved by Turing with an extremely general argument, a remarkably uniform idea that applies not only to his machine concept, but which applies generally with nearly every sufficiently robust concept of computability.''

\medskip
\end{enumerate}
The question we should like to consider is the extent to which these attributions to Turing are accurate. 

\section{The prima facie case against the Turing attribution}

The prima facie case against the Turing attribution, to be sure, consists of the observation that nearly all of the things attributed to Turing in the quotes above are not actually to be found in Turing's paper. He doesn't define or even discuss the halting problem as a decision problem; the phrase ``halting problem'' does not occur in his paper; there is no theorem in the paper called the Halting Theorem or any theorem or statement making a similar assertion; indeed, the word ``halt'' is absent; he does not discuss the halting of his machines at all, and makes no provision for the computational processes undertaken by his machines ever to stop; in particular, he has no convention as in contemporary accounts of a \emph{halt} state for the machines; none of the notation $\varphi_e(x){\downarrow}$, $\mathcal{K}_0$, and $\mathcal{W}_e$ occurs in Turing's paper, nor does any equivalent notation appear for these ideas; he doesn't use the undecidability of the halting problem to resolve the Entscheidungsproblem, but rather another undecidable decision problem; there are no remarks about the self-contemplative nature of Turing machines; and there is nothing like the self-referential proof of undecidability that we gave earlier to be found in Turing's paper. All his undecidability arguments proceed instead in multi-step reductions ultimately through the undecidability of his circle-free problem, which is not even computably equivalent to the halting problem but rather is strictly harder in the hierarchy of computational strength.

The Turing attribution for the undecidability of the halting problem has been challenged by a number of Turing scholars, including Jack Copeland \cite[p. 40]{copeland2004essential}, who explicitly claims that crediting Turing with stating and proving the halting theorem is erroneous. Copeland---as well as many others, e.g. Petzold \cite[p. 179]{petzold2008annotated}---views Martin Davis as the mathematician to whom we should attribute the result under discussion. This is primarily owed to Davis's influential book from 1958, \textit{Computability and Unsolvability}, where the phrase \textit{the halting problem} first appears in the literature:

\begin{quote}\small

[L]et $Z$ be a simple Turing machine. We may associate with $Z$ the following decision problem: 

\textit{To determine, of a given instantaneous description $\alpha$, whether or not there exists a computation of $Z$ that begins with $\alpha$.}

That is, we wish to determine whether or not $Z$, if placed in a given initial state, will eventually halt. We call this problem the \textit{halting problem} for $Z$. \cite[p. 70]{davis1958computability}
\end{quote}
This quote of Davis is followed by a proof of the undecidability of the halting problem for Turing machines. 

Notice, however, that Davis does not actually define here the halting problem as we would generally understand it today, that is, as a fully uniform decision problem. Rather, he has introduced separate individual halting problems for each particular machine $Z$. The uniform version, in contrast, would take input $\<M,x>$ and seek to determine whether machine $M$ halts on input $x$. Although the distinction between uniform and non-uniform computability is the crux of many subtle issues in computability, in this case the fully uniform halting problem is computably equivalent to the individual halting problem for a fixed machine $Z$, if it is a universal computer. If $Z$ on input $\<e,x>$ gives the same output behavior as program $e$ on input $x$, after all, then $e$ halts on $x$ if and only if $Z$ halts on $\<e,x>$. Another point of consideration is that if the halting problem for any individual machine $Z$ is undecidable, as Davis proves, then it follows immediately that the uniform version also is undecidable. 

It appears in any event that Davis wasn't actually the first to consider the halting problem, however, since Kleene already has an account of it in his classic book from 1952, \textit{Introduction to Metamathematics}, where he gives the familiar self-referential argument:
\begin{quote}\small
there is no algorithm for deciding whether any given machine, when started from any given initial situation, eventually stops. For if there were, then, given any number $x$, we could first decide whether $x$ is the Gödel number of machine $\mathfrak{M}_x$, and if so whether $\mathfrak{M}_x$ started scanning $x$ in standard position with the tape elsewhere blank eventually stops, and if so finally whether $x,\!1$ is scanned in standard position in the terminal situation. 
    \cite[p. 382]{kleene1952metamathematics}
\end{quote}
Thus, contrary to Copeland and Petzold, it would seem that Kleene has priority over Davis, although it was Davis who first used the words, ``halting problem.'' 

In his related historical analysis of the origins of the halting problem, Salvador Lucas comes to a similar conclusion, while also highlighting Alonzo Church's contributions:

\begin{quote}\small
    Martin Davis deserves credit for the wording of the current formulation of the halting problem (and also for the somehow standard expression \enquote{halting problem}). It seems, however, that it was Stephen C. Kleene the one who first mentioned the problem, although using a different wording. Alonzo Church made an important, early contribution by requiring termination as part of his notion of computation (\enquote{effective calculation}) in the $\lambda$-calculus. \cite[p. 8]{lucas2021origins}
\end{quote}

At any rate, in our nuanced final conclusion we shall explain why we find it perfectly reasonable to offer a qualified attribution for the undecidability of the halting problem to Turing.

\section{The circle-free problem}

The central undecidability result of Turing's paper concerns the problem of determining whether a given program is what he calls \emph{circle-free}, a concept arising in connection with his analysis of the computable real numbers. So let us go through his argument.\goodbreak

\subsection{Computable real numbers}

In the first sentence of his paper, Turing defines the computable real numbers as those for which there is a computable procedure to enumerate the decimal digits. He clarifies in \S2 that he intends to refer to binary notation, representing real numbers with their integer part and the sequence of their ``binary decimal'' digits, a sequence of $0$s and $1$s after the decimal point.

The first thing to say about this definition is that it is not the standard definition in use today in the subject of computable analysis, one of the numerous subjects to which Turing's paper gave rise. Turing's definition is problematic for several specific but subtle reasons, and adopting it would ultimately contradict several claims that Turing seems to make about the computable numbers.  

One major issue, for example, is that with Turing's definition, we cannot compute sums $a+b$ in a uniform computable manner---the mathematical fact of the matter is that there is no computable procedure to produce a program for enumerating the digits of $a+b$, given programs that enumerate the digits of $a$ and $b$, respectively, as input. Turing himself recognizes the issue in his corrective note \cite{Turing1938:On-computable-numbers-a-correction}, published a few years later. To see the problem, suppose we are performing the following sum, which for the purpose of this illustration we shall undertake with the familiar decimal digits, although the same phenomenon occurs in binary:
\begin{align*}
        &0.22222\ldots \\[-2pt]
 +\quad &0.77777\ldots \\[-10pt]
   \cline{1-2}
        &0.99999\ldots
\end{align*}
For the sum, it would be wrong for us automatically to begin enumerating the digits as $0.99999\ldots$ as shown, since we don't actually know that the input summands will continue in that pattern, and there could eventually be a carry term, if later digits add to ten or more, which would cause all these $9$s in the answer to roll over to $0$, with a $1$ in the front. That is, in this case the correct answer would have to begin $1.0000\ldots$, and so forth. But it would likewise also be wrong for us automatically to start with $1.000\ldots$, since the later digits of the summands might add up to less than $9$ in some place, making the answer strictly less than $1$, and so we should have started with the $9$s. The critical point is that we just can't determine even a single digit of this sum until we know whether there will be a carry or not, and this is a problem that is undecidable, computably equivalent to the halting problem. Using the Kleene recursion theorem, one can prove that there is no computable resolution of this issue (see \cite{Hamkins.blog2018:Alan-Turing-on-computable-numbers}, \cite[\S6.6]{Hamkins2021:Lectures-on-the-philosophy-of-mathematics}).\footnote{Namely, let $a$ be the program that enumerates $0.222\ldots$, with all $2$s, and let $b$ be the program that starts enumerating $0.777\ldots$, but at the same time runs the supposed program for $a+b$ until it sees what the first digits of the sum are. If that sum begins $0.999\ldots$, then $b$ should suddenly switch to all $8$s, which would make the answer wrong since there would be a carry; and if the sum begins $1.000\ldots$, then $b$ should suddenly switch to digit $0$s, which again would make the answer wrong. The Kleene recursion theorem is used to know that indeed there is such a self-referential program $b$ that can look at the sum of $a+b$ while it is computing.}

For similar reasons, we cannot computably multiply or apply any of the other standard analytic functions such as $e^x$, $\sin x$, $\ln x$, and so forth, which Turing says that he wants to do. The main problem is that although these are continuous functions, nevertheless the digits of the output reals for these functions, including addition, are not continuous in the digits of the input, because it can happen that very small changes in the input can cause a carry digit in the output, which will propagate to many digits prior.

The standard definition of computable real number used today in computable analysis, in contrast, is that a computable real number is one for which there is a computable sequence of rational numbers converging to it within a known computable rate of convergence---the $n$th approximation should be within $1/2^n$ of the limit. With this definition, all the problematic features of Turing's definition on this point are resolved. One can computably add computable real numbers, multiply them, compute trigonometric functions, exponentials and logarithms, in each case computing the result uniformly from the programs that compute the inputs to these functions.

Since the sequence of digits of a number does provide a convergent sequence of approximations, however, we may view Turing's original conception of computable real number simply as a slightly-too-strong notion. Nevertheless, his notion is sufficient to establish the undecidability of the various decision problems he considers, and so we shall proceed with it as he does.

Meanwhile, there is an interesting philosophical point to make in regard to this issue about extension versus intension, since the particular real numbers that are realized as computable under Turing's original definition and under the modified definition in computable analysis are exactly the same---the differences have not to do with the extension of the set of computable real numbers, but rather entirely with how those numbers are to be represented as programs for the purpose of further computing with them as input and output.

\subsection{Circle-freeness}

Turing defines that a program is \emph{circular} when it produces only finitely many binary output digits, and otherwise it is called \emph{circle-free}. So the computable real numbers are those computed by the circle-free programs, which produce the desired infinitely many binary digits on the write-only output tape.\footnote{He doesn't actually have a separate output tape, but models this by allowing the program on its one tape to mark certain cells as special, making them part of the output, and these will not be overwritten. In effect, it is like allowing a write-only output tape, but taking the output to consist only of the digits $0$ and $1$ on the output tape, ignoring other symbols or blank cells.}

Let us remark on the terminology. Although it is true that a program caught in an endlessly repeating computational loop will be able to produce only finitely many digits on the output tape, nevertheless this is not the only way for this to happen.  Consider, for example, an algorithm that plans to produce another output digit for each additional twin-prime instance that it finds, and so it begins searching for larger and larger twin-prime instances, producing another output digit for each new pair it finds. If there are infinitely many twin primes (an open question), then this program is circle-free. But if there are only finitely many twin primes, then the program will produce only finitely many digits. The issue we want to highlight is that in this second case, the program will count as ``circular'' in Turing's terminology, even though there is no sense in which the algorithm is caught in a repeating computational loop---rather, the program is computing furiously, searching hopelessly amongst larger and larger numbers for the sought next twin-prime instance without success. 

So, it is easy to object to the ``circular'' terminology on the grounds that it wrongly suggests that a program might fail to produce infinitely many digits only because it is caught in a repeating computational loop. So it will be best for us to keep in mind exactly the meaning that Turing has provided for these terms, namely, a program is circular, when it produces only finitely many digits of the output digit sequence, and circle-free, when it has succeeded in giving us an infinite digit sequence for the output real number. 

A rather worse problem is that we have observed people sometimes have a serious misunderstanding of Turing's terminology by taking ``circular'' on a superficial reading to mean that the program has failed to halt, which would lead to a misimpression that the circle-free problem is exactly the halting problem. 

It may be interesting to mention also that the repeating-computational-loop instances of nonhalting are actually computationally easy---semidecidable---since a repeating loop can be recognized for what it is when it occurs. Consequently, the difficulty of the halting problem lies rather entirely with the nonrepeating nonhalting instances. We can eventually recognize instances of halting and instances of repeating computational loops. What we cannot recognize are instances of nonhalting non-loopy computation. 

\subsection{The circle-free problem is computably undecidable}

Turing's main undecidability result is to show that there is no computable procedure to determine whether a given program is circle-free---in short, that what we may call the circle-free problem is undecidable. Let us briefly sketch the argument, which amounts to a computable analogue of Cantor's proof that the space $2^\N$ of all binary sequences is uncountable. Namely, we assume toward contradiction that the circle-free problem is computably decidable; we then use this fact to make a computable listing of all computable infinite binary sequences; and then, by flipping digits on the diagonal, we produce a computable binary sequence that isn't on the list, a contradiction. 

To begin, we assume toward contradiction that there is a computable procedure to decide if a given program is circle-free. Using this, we can produce 
a computable listing of all and only the circle-free programs $p_1$, $p_2$, $p_3$, and so forth. To do so, simply consider all the possible programs in turn, ordered by what Turing calls the ``description numbers'' of the programs, which are numbers encoding the full instruction details of the program, and then test each one, keeping only those that are circle-free on the list, the ``satisfactory'' programs, as Turing calls them. Now we are ready to define the diagonal real number $\beta$, whose binary representation is
  $$\beta=0.d_1d_2d_3\ldots$$
where the $n$th binary digit $d_n$ is obtained by flipping the $n$th binary digit of the binary sequence produced by $p_n$, that is, swapping $0$ for $1$ and conversely. The binary digit sequence of $\beta$ is computable, by the process we have just  described, but it cannot be the binary digit stream producing by any particular program $p_n$, since that digit stream differs from $\beta$ at the $n$th place $d_n$, since we had flipped that bit. This is a contradiction, and so the the circle-free decision problem must be computably undecidable.\footnote{Note that since we are diagonalizing against the binary sequences, rather than the real-number values of those sequences, there is no issue here about nonunique binary representations, such as with $0.111\ldots=1.000\ldots$, which occurs in some presentations of Cantor's theorem. In this sense, the proof is closer to Cantor's proof that $2^\N$ is uncountable than it is to the proof that $\R$ is uncountable.}

Following Peter Millican \cite{millican2021turing}, we remark that Turing's diagonal proof of the foregoing fact occurs in §8 of his paper, which starts with a paragraph that includes a subtle footnote to E. W. Hobson's \cite{hobson1921theory} book, referencing pages 87 and 88. Interestingly, these two pages contain an unexpected discussion of \textit{paradoxes of natural language definability}, such as those advanced by Julius König and Jules Richard. Based on this fact, plus the overall structure of Turing's paper and line of argumentation, Millican makes a forceful case that Turing likely stumbled upon the undecidability of the circle-free problem by experimenting with variations of Richard's paradoxical diagonal argument, in which the problematic notion of an \enquote{English-definable number} is replaced by Turing with the rigorous notion of a \enquote{computable number.}\footnote{Richard's paradox \enquote{shows} that the set of English-definable real numbers is both countable and uncountable. An example of such a number is $0.5555...$, for there is at least one English phrase which specifies it, e.g. \enquote{Point five recurring.} The overall argument starts with an observation that the set of finite English phrases which succeed in specifying a real number is countable (for they can be sorted by size and lexicographic order), hence the set of definable real numbers is countable as well. However, one can diagonally exhibit a real number which is not part of any proposed list that enumerates this set using the standard Cantorian strategy. Trivially, this generated number happens to be defined by the complex (yet finite) English phrase which merely describes the details of its diagonal construction. Putting these observations together entraps us into a contradiction. See also the first author's remarks on the related \emph{Math Tea Argument} in \cite{HamkinsLinetskyReitz2013:PointwiseDefinableModelsOfSetTheory}.} Given Turing's mathematically precise understanding of computable numbers in terms of his proposed model of computation, a Richard-style diagonal argument for the case of computable numbers can only work now if there would be a computable procedure to determine whether a given program is circle-free. 
Hence, on the pain of paradox, such a program cannot exist.\goodbreak

Turing gives a second proof of the computable undecidability of the circle-free problem in the latter part of \cite[\S8]{Turing1936:On-computable-numbers}, proceeding like this. If circle-freeness were decidable, then we could design a program that printed the first output digit of the first circle-free program, then the second output digit of the second circle-free program, and continuing with the $n$th output digit of the $n$th circle-free program for every $n$. This program would be circle-free, and so it would be the $n$th circle-free program for some $n$. But then at stage $n$ of the process, the algorithm would be left waiting to find out what its own $n$th output digit is before saying what it is. Under the assumption that circle-freeness is computability decidable, Turing has thus produced a program that is circle-free and yet cannot be. Contradiction.

\subsection{Circle-free problem is harder than the halting problem}

We would like to remark on a rather curious aspect about the situation here. Although as we mentioned, Turing is commonly credited with proving the undecidability of the halting problem, nevertheless the central undecidability result in his paper concerns the circle-free problem, which in fact is strictly harder than the halting problem in the hierarchy of Turing degrees.\footnote{We thank an anonymous referee for pointing out that this was first conjectured by \cite{Post1947:Recursive-Unsolvability-of-a-problem-of-Thue}. A proof of this result can also be found in \cite[§5.2]{lucas2021origins}.} The claim that it is undecidable, therefore, would be a strictly weaker result.

Let us explain. The circle-free decision problem has a natural logical complexity of $\Pi^0_2$, that is, with quantifier complexity $\forall\exists$, since a given program is circle-free if and only if for every natural number $n$, there is a stage at which the program has produced at least $n$ output digits. The halting problem, in contrast, has simpler complexity $\Sigma^0_1$, an arithmetical existential $\exists$ assertion, since any given instance of halting is witnessed by the length of the halting computation itself.

Indeed, the circle-free problem is a complete $\Pi^0_2$ problem, meaning that every $\forall\exists$ problem reduces to the circle-free problem, and this shows that the $\Pi^0_2$ classification cannot be simplified. To see this, let us reduce an arbitrary $\Pi^0_2$ assertion $\forall n\exists k\ \varphi(n,k,x)$, where $\varphi$ has only bounded quantifiers, to the circle-free problem. Let $e$ be the program which on input $x$ systematically considers $n=0$ and then $n=1$ and $n=2$ and so on in turn. For each $n$, it looks for a $k$ for which $\varphi(n,k,x)$. If found, then $e$ produces another digit on the output tape, say, digit $1$, and then moves on to $n+1$. The original assertion $\forall n\exists k\ \varphi(n,k,x)$ is true for $x$ if and only if program $e$ is circle-free. So we have reduced any given $\Pi^0_2$ statement to the circle-free problem. 

In particular, this means that the circle-free problem is Turing equivalent not to the halting problem, but to the double jump $0''$, the double halting problem, that is, the halting problem relativized to the halting problem. 

\subsection{Refining the argument}

Nevertheless, an easy modification of the diagonalization argument relies only on a weaker version of the circle-free problem, namely, the question whether a given program will produce at least a given finite number of digits on the output tape. Let us call this \emph{circle-free for $n$ digits}. 

Turing could have undertaken his diagonalization argument using only this weaker decision problem. Namely, with a solution to this problem, we could enumerate all possible programs, and cut down to a list of those programs $p_1$, $p_2$, $p_3$, and so forth, such that $p_n$ is circle-free for at least $n$ digits. The idea is that this list would include all the fully circle-free programs, which are the ones we want, but it may also include 
some other programs, which may produce only finitely many digits, but enough to cover the diagonal digit used in the diagonalization procedure of the argument. Now, with this list we define the diagonal binary sequence $\beta$ as before, flipping the $n$th digit in the sequence produced by program $p_n$, and so this will be a computable infinite binary sequence that is different from every computable infinite binary sequence, a contradiction. 

Our point is that this weaker problem, which would have sufficed in Turing's diagonal argument, is computably equivalent to the halting problem. The question whether a given program $e$ produces $n$ digits on the output tape, for example, is reducible to the halting problem, since for any $e$ and $n$ we could design a program $p$ that runs $e$ until $n$ digits appear and then halts. Program $p$ halts if and only if $e$ produces at least $n$ digits. Conversely, the halting problem is reducible to the circle-free for $n$ digits problem, since given any program $e$ to be run on input $x$, we can design a program $q$ that runs $e$ on $x$ and if this halts, afterwards produces the digits of $\pi$. So $e$ halts on input $x$ if and only if program $q$ produces at least $1$ digit, or at least $n$ digits for any particular $n$.

Alternatively, Turing could have used instead the \emph{one-more-digit} problem, which asks of a given program $p$ and time number $t$, whether $p$ will produce another digit after time $t$. Having a computable method to solve this problem would be enough to perform the Turing diagonalization, since at stage $n$ we can test the next candidate program to see if it will produce $n$ digits by running it and repeatedly asking whether it will produce one more digit, until either we have $n$ digits or we know that it will not produce $n$ digits (in which case that program drops off the list and we move on to the next program). In this way, Turing's argument shows that the one-more-digit problem also is computably undecidable.

And again the point is that this one-more-digit problem, which would have sufficed in Turing's undecidability argument, is computably equivalent to the halting problem. It reduces to the halting problem, since $e$ will produce at least one more digit after time $t$ if and only if the program $p$ that we would design to halt in that event actually halts. And conversely, the halting problem reduces to the one-more-digit problem, since $e$ halts on input $x$ if and only if the $\pi$ program $q$ we considered above produces one more digit after time $0$.

Turing doesn't use either of these problems, however, but takes instead a different approach.

\section{The printing problem}\label{Section.printing}

Turing uses a very clever method to show a certain other problem is undecidable, what we shall call the \emph{printing} problem, and this problem is easily seen to be computably equivalent to the halting problem. Specifically, Turing proves there can be no computable procedure to determine whether a given program will ever print a certain symbol as output, say, the symbol $0$. The problem is very close to the halting problem, since printing a certain special symbol can be taken as a triggering event for halting, and so we can immediately reduce the printing problem to the halting problem. And conversely, we can reduce halting to printing, since for any given program/input combination, we design a new program that would operate exactly the same, except that if a halting event should occur, it makes beforehand a one-time printing of a special character right at that moment. In this way, Turing's printing problem and the halting problem can be seen as close natural variations of one another, almost identical.

\subsection{Alternative halting criterion}

\newcommand\halt{\boldsymbol{\downarrow}}
To press this point further, let us imagine a computationally equivalent model of Turing machines in which halting is determined not by reaching a special halting state, as usual in contemporary accounts, but rather by the printing of a special \emph{halt} symbol, perhaps the symbol $\halt$\,. That is, for this model of computability, a computation would stop, by definition, exactly in the event that $\halt$ is printed. For this model, the halting problem becomes identical to the printing problem for this special symbol. In our view, this convention would in many respects appear to be a more natural approach to halting than the standard halt-state convention, and in any case, as we have mentioned, the halt-state convention does not appear in \cite{Turing1936:On-computable-numbers}---being focused on the computation of infinite binary sequences, he simply does not specify any criterion at all for the halting of a computation. In this sense, one may regard the printing problem to be a very close variant of the halting problem. Indeed, with only the slight variation in machine halting formalism that we have suggested, the entire field of computability theory might have been talking everywhere about the printing problem instead of the halting problem, since these two problems can in principle play the same fundamental role in the theory.

\subsection{Turing's clever proof of undecidability}

Let us now explain Turing's clever argument proving the computable undecidability of the printing problem. He mounts an unusual kind of reduction, showing that if the printing problem were decidable, then also the circle-free problem would be decidable, which he had already proved is not the case. 

What is both interesting and unusual about this argument is that it is not a simple relative computability result as might be expected, showing that the circle-free problem is computable relative to the printing problem, since in fact no such reduction is possible. The circle-free problem, after all, has complexity complete-$\Pi^0_2$, which cannot be reduced to a $\Sigma^0_1$-problem such as the printing problem. Nevertheless, Turing does argue that if we could actually solve the printing problem computably, then we would be able computably to solve the circle-free problem, which we cannot. And so, he concludes, the printing problem is undecidable. 

He argues as follows. Assume we had a computable procedure for deciding whether a given program would ever print the symbol 0 on the (write-only) output tape. Next, we show on this basis that the \emph{infinite-symbol} problem is computably undecidable, the problem whether a given program will produce infinitely many of the given symbol on the output tape (this is a $\Pi^0_2$-complete problem). To see this, for any given program $p$, we can systematically produce programs $p_k$ that would replace the first $k$ instances of $0$ on the output tape, if these occur, with a different symbol, say, $\bar 0$. Now, we consider a program $q$ that  would consider each $p_k$ in turn, asking whether it ever prints $0$. If yes, we move on to the next $p_{k+1}$. If no, however, if some $p_k$ is found that would never print $0$, then $q$ will print 0 at that moment (and only in this kind of case). Finally, we don't actually run the program $q$, but instead observe that $p$ prints infinitely many 0s just in case $q$ does not print a $0$, since the only way $q$ prints a zero is if some $p_k$ does not print a $0$, which would mean that the original $p$ printed at most $k$ zeros. 

This is not a straightforward reduction of one problem to another, but rather an argument that if one problem were actually computably decidable, then so would be the other. At bottom, the logic of the argument is like this: if we had a computable way of finding whether existential statements are true, then we could iterate this with negation to also compute $\forall\exists$ assertions, since $\forall k\exists n\,\varphi$ fails just in case there is some $k$ for which the existential statement about it fails. In short, if in general existential statements are decidable, then the whole arithmetic hierarchy collapses. 

Thus, Turing has proved that if the printing problem were computably decidable, then so also is the infinite-symbol problem. But now the circle-free problem is easily computable from this, since a program produces infinitely many output binary digits just in case it produces either infinitely many $0$s or infinitely many $1$s. This would contradict the earlier proof that the circle-free problem is not decidable, and so, therefore, the printing problem is not computably decidable. 

In light of the very close connection between the printing problem and the halting problem, the undecidability of the printing problem is perhaps the closest Turing gets in his paper to proving the undecidability of the halting problem. And indeed, in our view, it is very close.

\subsection{A simple self-referential proof of undecidability for printing}

Turing thus showed that the printing problem is undecidable by mounting a reduction to and through the undecidability of the circle-free problem. But let us illustrate how one may improve upon Turing with a simpler self-referential proof of the undecidability of the printing problem in the style of the standard contemporary proof of the undecidability of the halting problem. There was actually no need for Turing's detour through the circle-free problem. 

Namely, assume toward contradiction that the printing problem were computably decidable, and fix a method of solving this problem. Using this as a subroutine, consider the algorithm $q$ which on input $p$, a program, asks whether $p$ on input $p$ would ever print $0$ as output. If so, then $q$ will halt immediately without printing $0$; but if not, then $q$ prints $0$ immediately as output. So $q$ has the opposite behavior on input $p$ with respect to printing $0$ as output than $p$ has on input $p$. Running $q$ on input $q$ will therefore print $0$ as output if and only if it will not, a contradiction. So the printing problem is computably undecidable. 

The fundamental similarity of this argument with that for the halting problem, in our view, buttresses our claim that these two decision problems are close variants, almost identical.

\section{The Entscheidungsproblem}

Turing used the undecidability of the printing problem to establish the undecidability of the Hilbert-Ackermann \cite{hilbert1931grundzuge} Entscheidungsproblem. 

Sometimes one hears an alternative variant of the Entscheidungsproblem, described as the problem of deciding whether a given arithmetic statement is true:

\begin{quote}\small
    [I]s there a definite procedure that can be applied to every statement that will tell us in finite time whether or not the statement is true or false? The idea here is that you could come up with a mathematical statement such as, \enquote{Every even integer greater than 2 can be expressed as the sum of two prime numbers,} hand it to a mathematician (or computer), who would apply a precise recipe (a \enquote{definite procedure}), which would yield the correct answer \enquote{true} or \enquote{false} in finite time. (...) [This] is known by its German name as the \textit{Entscheidungsproblem} (...). \cite[pp. 58-59]{mitchell2009complexity}
\end{quote} 
Similarly:
\begin{quote}\small
    The \emph{Entscheidungsproblem} was the dream, enunciated by David Hilbert in the 1920s, of designing a mechanical procedure to determine the truth or falsehood of any well-formed mathematical
statement. \cite[\S3.1]{Aaronson2013:Why-philosopher-should-care-about-computational-complexity}
\end{quote}

To show that this version of the problem is undecidable, it would suffice to use the circle-free problem directly, since every instance of circle-freeness is an instance of whether a certain $\forall\exists$ arithmetic assertion is true. We can therefore have no computable algorithm to decide whether such assertions are true, since there is no computable algorithm to decide whether a given program is circle-free.\footnote{On this note, Mitchell claims in the same book chapter that Turing showed that the Entscheidungsproblem is undecidable by considering the halting problem as the relevant mathematical statement: \enquote{Thus---and this is Turing's big result---\textbf{there can be no definite procedure for solving the Halting problem}. The Halting problem is an example that proves that the answer to the \textit{Entscheidungsproblem} is \enquote{no}; not every mathematical statement has a definite procedure that can decide its truth or falsity.} (p. 68, bold original)}

But in our view it is probably more correct historically to understand the Entscheidungsproblem as the problem of deciding whether a given sentence is a logical consequence of a given theory (as specified by a computable list of axioms), that is, the validity problem. By the completeness theorem of first-order logic, a statement is a logical consequence of a theory if and only if it is provable in that theory, and so the validity problem is the same as the provability problem, the problem of deciding whether a given sentence is provable from a given (computable) theory. 

Turing's doctoral supervisor, Alonzo Church, describes the problem as follows:\footnote{Also quoted in \cite[p. 53]{copeland2004essential}, where one can also find a rich discussion of neighboring issues pertaining to the Entscheidungsproblem.}

\begin{quote}\small

By the Entscheidungsproblem of a system of symbolic logic is here understood the problem to find an effective method by which, given any expression $Q$ in the notation of the system, it can be determined whether or not $Q$ is provable in the system. \cite[footnote 6]{church1936note}
    
\end{quote}

Turing construes the Entscheidungsproblem problem as follows, matching the Hilbert-Ackermann variant:

\begin{quote}\small

\hspace{1em} The results of §8 have some important applications. In particular, they can be used to show that the Hilbert Entscheidungsproblem can have no solution. For the formulation of this problem I must refer the reader to Hilbert and Ackermann's \textit{Grundzüge der Theoretischen Logik} (Berlin, 1931), chapter 3. 

I propose, therefore, to show that there can be no general process for determining whether a given formula $\mathfrak{U}$ of the functional calculus $\mathsf{K}$ is provable, \textit{i.e.} that there can be no machine which, supplied with any one $\mathfrak{U}$ of these formulae, will eventually say whether $\mathfrak{U}$ is provable. \cite[§11, p. 259]{Turing1936:On-computable-numbers}
\end{quote}

Tangentially, the phrase \textit{functional calculus $\mathsf{K}$} is seldom used in the contemporary literature, but it is essentially synonymous with classical first-order logic (without identity). Curiously, this non-standard terminology isn't even there in the edition of the book that Turing cites---Ackermann himself says in the preface of the book that the technical vocabulary was updated to reflect modern usage.

To prove that the Entscheidungsproblem is not computably decidable, Turing reduces the printing problem to it. This method would not work directly with the circle-free problem, since a $\Pi^0_2$ problem like circle-freeness simply cannot computably reduce to a $\Sigma^0_1$ problem like validity (although one can imagine the possibility of another instance here of Turing's clever method mentioned in \S\ref{Section.printing}).

In contemporary language, Turing's argument amounts to a version of the fact that Peano arithmetic PA proves all true existential statements, and so let us cast it in this way. If a given program does eventually print symbol $0$ at some point, then this will be provable in PA, since the computational history to that juncture would be finite and we could take the entire computational history down as a big sequence and prove that indeed the computation proceeded exactly as indicated according to it and thus led ultimately to a $0$ being printed. Conversely, if a program does not ever print $0$, then this will be consistent with PA, since this is true in the standard model $\N$, and so PA will not prove that it prints a $0$. In short, a program will print a $0$ if and only if PA proves that it will. So provability in PA cannot be computably decidable, since then the printing problem would be as well, which Turing proved it is not.

\section{Mathematical attribution practice}

Let us briefly discuss a certain cultural observation about the practice of attributions in mathematics, namely, that mathematicians often exhibit a certain generosity in their attributions. What we mean is that mathematicians often attribute results, ideas, and methods to earlier thinkers in instances where those ideas, results, and methods do not actually appear in the original cited work, provided that those ideas, results, and methods grew directly out of that work. Thus, a mathematician will sometimes be credited for unstated corollaries of their work or even for theorems, which they did not state, but which can be easily proved using the methods and ideas that indeed they had provided.

Such, we claim, is the situation of Alan Turing and the halting problem. He hadn't stated the halting problem or proved that it is undecidable, but he provided all the tools and ideas that we needed to do so ourselves, and far more. It is a small step, after all, from the undecidability of printing to the undecidability of the halting problem.

Let us exhibit a few examples of this cultural practice.

\begin{enumerate}
\item We commonly attribute the irrationality of $\sqrt{2}$ to the Pythagoreans, although they would have had a far more geometric understanding of the result. The most common contemporary proof, a short sequence of number-theoretic algebraic equations and parity considerations, uses notation and concepts that did not arise until many centuries after the Pythagoreans.

\item The Chinese remainder theorem is often credited to 5th century mathematician Sunzi, although that early work consisted merely of examples, without a general result, which was not formulated until many centuries later. 

\item Euler is commonly credited with inventing graph theory, particularly in his 1736 solution to the K\"onigsberg bridge problem (see \cite{BiggsLloydWilson1986:Graph-theory-1736-1936}), but his paper does not present graph theory with our current abstract perspective on the subject. 

\item The Cayley-Hamilton theorem in linear algebra, asserting that every square matrix over a commutative ring fulfills its characteristic equation, was originally proved by Hamilton 1853 only for quaternions, and stated by Cayley 1858 for the $3\times 3$ real case, proved only for the $2\times 2$ case. Cayley \cite{Cayley1858:A-memoir-on-the-theory-of-matrices} writes “..., I have not thought it necessary to undertake the labor of a formal proof of the theorem in the general case of a matrix of any degree,” although the general case was indeed proved in Frobenius in 1878. 

\item We commonly attribute the fundamental theorem of finite games, also known as ``Zermelo's theorem,'' to Zermelo 1913 (see \cite{Larson2010:Zermelo-1913,Hamkins2025:The-fundamental-theorem-of-finite-games}), the result stating that in any finite two-player game of perfect information, either one of the players has a winning strategy or both players have drawing strategies. Zermelo does not provide the full contemporary understanding of what constitutes a finite game or a two-player game of perfect information, being rather concerned with board games like chess and the fact that there are only finitely many positions on the board. Zermelo is not commonly attributed with the open-determinacy result of Gale and Stewart \cite{GaleStewart1953:InfiniteGamesWithPerfectInformation}, even though several standard proofs of open determinacy are no more difficult than and indeed almost identical to proofs of the finite-games result---Zermelo simply wasn't working with infinite games. 

\item We commonly attribute to \Godel\ various strong statements of the first incompleteness theorem, even though the theorem actually appearing in his work involves extra hypotheses such as $\omega$-consistency, which are now seen as extraneous.

\end{enumerate}

In the case of Turing, we have already mentioned Odifreddi's \cite[p. 150]{odifreddi1992classical} statement of the undecidability of the halting problem, as follows: 

\bigskip 

\noindent\textbf{Theorem II.2.7 Unsolvability of the Halting Problem (Turing [1936])} \textit{The set defined by $\braket{x,e} \in \mathcal{K}_0 \iff x \in \mathcal{W}_e \iff \varphi_e(x)\downarrow$ is r.e. and nonrecursive.}

\bigskip 

In this case, Turing is explicitly credited with this result, even though none of the $\mathcal{K}_0, \varphi_e(x)\downarrow,\mathcal{W}_e$ notations, or indeed even the concept of r.e.~or the enumeration of the r.e.~sets were part of the 1936 paper. 

We refer the reader also to \cite{CallardYourgrau2025:Who-Proved-the-Independence-of-the-ContinuumHypothesis} for the concept of \emph{credit nihilism}, according to which the general practice of assigning credit for ideas is intellectually and even morally bankrupt. 

\section{A possible explanation for the misattributions}\label{Section.Possible-explanation}

It is easy enough to see how a student or beginning researcher might fall into a mistaken impression that Turing had proved the undecidability of the halting problem. Looking at the paper, one finds that he discusses at length the concept of \emph{circular} machines and \emph{circle-free} machines, and with a light reading of what is in truth a technically demanding paper full of finicky details it may seem natural to slip into a misunderstanding of this terminology by which ``circular'' would mean that the computational process is caught in a repeating loop or otherwise not giving output, hence not halting; and so ``circle-free'' would therefore mean the opposite: halting. By this reading, therefore, the circle-free problem would be misinterpreted as identical to the halting problem. And since Turing definitely proves the undecidability of this circle-free problem, someone with this misunderstanding of the terminology would take him to have proved the undecidability of the halting problem. 

We emphasize that this would indeed be a severe misunderstanding of Turing's terminology and the content of his paper, since circle-free for Turing means that the machine has succeeded in producing an infinite binary output stream, which makes it a $\Pi^0_2$-complete decision problem, strictly harder than the halting problem. The circle-free problem is simply not the same as nor even computably equivalent to the halting problem.

Perhaps this particular misunderstanding of Turing's circle-free terminology, however, may be the source of some of the misattributions? Although we have definitely observed this mistake amongst students and in discussions online, it wouldn't seem to occur as  readily amongst expert researchers, such as those authors making the mistaken Turing attributions in the quotations we mentioned in \S\ref{Section.Quotes-giving-Turing-attribution}.

\section{A nuanced conclusion}

In conclusion, let us offer a nuanced account. Strictly speaking, Turing did not prove nor even state the undecidability of the halting problem in his 1936 paper \cite{Turing1936:On-computable-numbers}, and it is incorrect to suggest that this result or any discussion of it can be found there. It is especially incorrect to attribute to Turing \cite{Turing1936:On-computable-numbers} the common self-referential proof of the undecidability of the halting problem, since nothing like that argument appears in Turing's paper. 

Nevertheless, Turing did provide a robust framework of ideas sufficient to lead directly to the undecidability of the halting problem. More than this, he proved the undecidability of the printing problem, which is easily seen to be computably equivalent to the halting problem, and indeed, perhaps a mere translation of it---one can merely view halting as a trigger event to print a special \emph{halt} symbol and conversely. An alternative model of Turing machine behavior could simply replace the halt-state manner of halting with a halt-printing criterion, which would make the halting problem an instance of the printing problem. In this way, the printing problem can be seen as a close variant of the halting problem.\goodbreak

In light of all these considerations, we would certainly find it correct to say:
\begin{quote}\it
Turing \cite{Turing1936:On-computable-numbers} provided all the core ideas leading to a proof of the undecidability of the halting problem.
\end{quote}
But more than this, we would find it completely fine to attribute the undecidability of the halting problem to Turing in the following only slightly qualified way:
\begin{quote}\it 
Turing essentially proved the undecidability of the halting problem in \cite{Turing1936:On-computable-numbers}.
\end{quote}
What he proved, of course, is the undecidability of the printing problem, but we have described how this can be regarded as a simple translation of the halting problem. There are essentially no substantive ideas required to transform Turing's proof of the undecidability of the printing problem to a proof of the undecidability of the halting problem. In this sense, he essentially proved the undecidabilty of the halting problem (even if he did not do so by the standard self-referential argument commonly used today).

Perhaps it is a kind of historical accident that so much of computability theory has become founded on the halting problem rather than on Turing's printing problem, since every use of the halting problem in computability theory to our knowledge can be easily be undertaken with transparent, superficial changes to use instead the printing problem. We might all have been talking about the printing problem everywhere instead of the halting problem, and the resulting computability theory would be essentially the same.

\printbibliography

\end{document}

%% file: MathMacrosJDH.tex
\newtheorem*{theorem*}{Theorem}

\newtheorem*{maintheorem*}{Main Theorem}
\newtheorem*{maintheorems*}{Main Theorems}

\newtheorem*{corollary*}{Corollary}
\newtheorem*{corollaries*}{Corollaries}

\theoremstyle{definition}

\newtheorem*{definition*}{Definition}

\newtheorem*{question*}{Question}

\newtheorem*{questions*}{Questions}
\newtheorem*{mainquestion*}{Main Question} 
\newtheorem*{openquestion*}{Open Question} 
\theoremstyle{remark}

\newcommand{\QED}{\end{proof}}

\def\proclaim[#1]{{\bf #1}}
\def\BF#1.{{\bf #1.}}

\def\says#1:#2\par{\item[#1] #2\par}

\newcommand{\Godel}{G\"odel}

\newcommand{\N}{{\mathbb N}}

\newcommand{\R}{{\mathbb R}}

\makeatletter
\newcommand{\dotminus}{\mathbin{\text{\@dotminus}}}
\newcommand{\@dotminus}{%
  \ooalign{\hidewidth\raise1ex\hbox{.}\hidewidth\cr$\m@th-$\cr}%
}
\makeatother
\newcommand\dbrace{\hskip-1.5em\raise3pt\hbox{\rotatebox[origin=c]{-35}{$\left.\strut^{\phantom{|}}\right\}$}}}

\newcommand\UParroW{{\setbox0\hbox{$\Uparrow$}\rlap{\hbox to \wd0{\hss$\mid$\hss}}\box0}}

\renewcommand{\setminus}{\raise.3ex\hbox{\rotatebox{-20}{$-$}}} 

\newcommand{\smalllt}{\mathrel{\mathchoice{\raise2pt\hbox{$\scriptstyle<$}}{\raise1pt\hbox{$\scriptstyle<$}}{\raise0pt\hbox{$\scriptscriptstyle<$}}{\scriptscriptstyle<}}}
\newcommand{\smallleq}{\mathrel{\mathchoice{\raise2pt\hbox{$\scriptstyle\leq$}}{\raise1pt\hbox{$\scriptstyle\leq$}}{\raise1pt\hbox{$\scriptscriptstyle\leq$}}{\scriptscriptstyle\leq}}}

   \def\DHLhksqrt#1#2{%
   \setbox0=\hbox{$#1\sqrt{#2\,}$}\dimen0=\ht0
   \advance\dimen0-0.2\ht0
   \setbox2=\hbox{\vrule height\ht0 depth -\dimen0}%
   {\box0\lower0.4pt\box2}}

\def\[#1]{\mathopen{\lbrack\!\lbrack}#1\mathclose{\rbrack\!\rbrack}}
\newbox\gnBoxA
\newbox\gnBoxB
\newdimen\gnCornerHgt
\setbox\gnBoxA=\hbox{\tiny$\ulcorner$}
\global\gnCornerHgt=\ht\gnBoxA
\newdimen\gnArgHgt
\def\gcode #1{%
\setbox\gnBoxA=\hbox{$#1$}%
\setbox\gnBoxB=\hbox{$\bar #1$}%
\gnArgHgt=\ht\gnBoxB%
\ifnum     \gnArgHgt<\gnCornerHgt \gnArgHgt=0pt%
\else \advance \gnArgHgt by -\gnCornerHgt%
\fi \raise\gnArgHgt\hbox{\tiny$\ulcorner$} \box\gnBoxA %
\raise\gnArgHgt\hbox{\tiny$\urcorner$}}
\newcommand{\UnderTilde}[1]{{\setbox1=\hbox{$#1$}\baselineskip=0pt\vtop{\hbox{$#1$}\hbox to\wd1{\hfil$\sim$\hfil}}}{}}
\newcommand{\Undertilde}[1]{{\setbox1=\hbox{$#1$}\baselineskip=0pt\vtop{\hbox{$#1$}\hbox to\wd1{\hfil$\scriptstyle\sim$\hfil}}}{}}
\newcommand{\undertilde}[1]{{\setbox1=\hbox{$#1$}\baselineskip=0pt\vtop{\hbox{$#1$}\hbox to\wd1{\hfil$\scriptscriptstyle\sim$\hfil}}}{}}
\newcommand{\UnderdTilde}[1]{{\setbox1=\hbox{$#1$}\baselineskip=0pt\vtop{\hbox{$#1$}\hbox to\wd1{\hfil$\approx$\hfil}}}{}}
\newcommand{\Underdtilde}[1]{{\setbox1=\hbox{$#1$}\baselineskip=0pt\vtop{\hbox{$#1$}\hbox to\wd1{\hfil\scriptsize$\approx$\hfil}}}{}}

\renewcommand{\iff}{\mathrel{\leftrightarrow}}

\def\<#1>{\left\langle#1\right\rangle}

\newcommand{\cell}[1]{\boxit{\hbox to 17pt{\strut\hfil$#1$\hfil}}}
\newcommand{\head}[2]{\lower2pt\vbox{\hbox{\strut\footnotesize\it\hskip3pt#2}\boxit{\cell#1}}}
\newcommand{\boxit}[1]{\setbox4=\hbox{\kern2pt#1\kern2pt}\hbox{\vrule\vbox{\hrule\kern2pt\box4\kern2pt\hrule}\vrule}}
\newcommand{\Col}[3]{\hbox{\vbox{\baselineskip=0pt\parskip=0pt\cell#1\cell#2\cell#3}}}
\newcommand{\tapenames}{\raise 5pt\vbox to .7in{\hbox to .8in{\it\hfill input: \strut}\vfill\hbox to
.8in{\it\hfill scratch: \strut}\vfill\hbox to .8in{\it\hfill output: \strut}}}
\newcommand{\Head}[4]{\lower2pt\vbox{\hbox to25pt{\strut\footnotesize\it\hfill#4\hfill}\boxit{\Col#1#2#3}}}
\newcommand{\Dots}{\raise 5pt\vbox to .7in{\hbox{\ $\cdots$\strut}\vfill\hbox{\ $\cdots$\strut}\vfill\hbox{\
$\cdots$\strut}}}